\documentclass[letterpaper, 10 pt, conference]{ieeeconf}  % Comment this line out if you need a4paper

\IEEEoverridecommandlockouts                              % This command is only needed if
\include{TCILATEX}

\usepackage{graphics} % for pdf, bitmapped graphics files
\usepackage{epsfig} % for postscript graphics files
\usepackage{amsmath} % assumes amsmath package installed
\usepackage{amssymb}  % assumes amsmath package installed
\usepackage{mathrsfs}
\usepackage{nicefrac}

\newcommand{\Sabs}[1]{\left\lceil #1 \right\rfloor}

\newcommand{\abs}[1]{\left| #1 \right|}
\newcommand{\norm}[1]{\left|\left| #1 \right|\right|}

\newcommand{\sign}{\textup{sign}}

\newtheorem{lemm}{Lemma}

\title{\LARGE \bf
Motion Planning and Robust Tracking for the Heat Equation using Boundary Control
}

\author{Diego Guti\'errez-Oribio$^{1}$, Yury Orlov$^{2}$, Ioannis Stefanou$^{1}$ and Franck Plestan$^{3}$% <-this % stops a space
\thanks{$^{1}$ Nantes Universit\'e, \'Ecole Centrale de Nantes, CNRS, GeM, UMR 6183, F-44000 Nantes, France
        ({\tt\small diego.gutierrez-oribio@ec-nantes.fr} and {\tt\small ioannis.stefanou@ec-nantes.fr})}%
\thanks{$^{2}$ CICESE Research Center, 22860 Ensenada, M\'exico
        ({\tt\small yorlov@cicese.mx})}%
\thanks{$^{3}$ Nantes Universit\'e, \'Ecole Centrale de Nantes, CNRS, LS2N, UMR 6004, F-44000 Nantes, France
        ({\tt\small franck.plestan@ec-nantes.fr})}%
}

\begin{document}

\maketitle
\thispagestyle{empty}
\pagestyle{empty}

%%%%%%%%%%%%%%%%%%%%%%%%%%%%%%%%%%%%%%%%%%%%%%%%%%%%%%%%%%%%%%%%%%%%%%%%%%%%%%%%
\begin{abstract}

Robust output tracking is addressed in this paper for a heat equation with Neumann boundary conditions and anti-collocated boundary input and output. The desired reference tracking is solved using the well-known flatness and Lyapunov approaches. The reference profile is obtained by solving the motion planning problem for the nominal plant. To robustify the closed-loop system in the presence of the disturbances and uncertainties, it is then augmented with PI feedback plus a discontinuous component responsible for rejecting matched disturbances with \textit{a priori} known magnitude bounds. Such control law only requires the information of the system at the same boundary as the control input is located. The resulting dynamic controller globally exponentially stabilizes the error dynamics while also attenuating the influence of Lipschitz-in-time external disturbances and parameter uncertainties. For the case when the motion planning is performed over the uncertain plant, an exponential Input-to-State Stability is obtained, preserving the boundedness of the tracking error norm. The proposed controller relies on a discontinuous term that however passes through an integrator, thereby minimizing the chattering effect in the plant dynamics. The performance of the closed-loop system, thus designed, is illustrated in simulations under different kinds of reference trajectories in the presence of external disturbances and parameter uncertainties.

\end{abstract}

%%%%%%%%%%%%%%%%%%%%%%%%%%%%%%%%%%%%%%%%%%%%%%%%%%%%%%%%%%%%%%%%%%%%%%%%%%%%%%%%
\section{INTRODUCTION}

The heat equation is a first-order parabolic partial differential equation used to describe the heat and  other diffusion processes such as chemical reactions, population growth, market price fluctuations, and fluids flow to name a few.

The motion planning has been addressed for the heat equation due to the flat output approach capable of parameterising the state dynamics using it as a reference \cite{b:Laroche-Martin-Rouchon-2000,b:Lynch-Rudolph-2002,b:Krstic-Smyshlyaev-2008,b:Krstic-Magnis-Vazquez-2009,b:Meurer-Kugi-2009}. Such approach provides an open-loop control able to track the desired reference but only for the nominal system with no uncertainties  and disturbances and starting at some specific initial condition.

In order to robustify such a tracking  control, one may invoke a feedback control. Examples of boundary tracking control of the heat equation utilize output feedback \cite{b:Jin-Guo-2018,b:Wu-Feng-2020} and  backstepping designs \cite{b:Krstic-Smyshlyaev-2008,b:Krstic-Magnis-Vazquez-2009,b:Meurer-Kugi-2009} as well as sliding-modes techniques \cite{b:Pisano-Orlov-Usai-2011,b:Pisano-Orlov-2012}.

Clearly, more realistic models should count for uncertainties and disturbances. A classic approach to deal with constant disturbances in the finite dimension setting is the integral action \cite{b:Khalil2002}. Dealing with a wider class of disturbances, sliding-mode control has long been recognized as a powerful control method to counteract non-vanishing external disturbances and unmodelled dynamics, even for infinite dimension systems (see, \textit{e.g.}, \cite{b:Orlov-Utkin-1987} and a very recent monograph \cite{b:Orlov-2020}).

The design of robust controllers in the heat equation has been addressed in \cite{b:Pisano-Orlov-2012} using sliding-modes and using $H_\infty$ control in \cite{b:Fridman-Orlov-2009}. For the case of the tracking task, in \cite{b:Pisano-Orlov-Usai-2011} a control has been designed to compensate Lipschitz-in-time disturbances using sliding-modes, but using distributed control, and in \cite{b:Jin-Guo-2018,b:Wu-Feng-2020} robust boundary controls have been designed but only to compensate bounded disturbances and requiring to design two external systems (disturbance estimator and servo system) to fulfil the task. Therefore, the design of boundary controllers capable of the boundary tracking of the heat equation under a wider kind of disturbances calls for further investigation.

In this paper, a simple control strategy is proposed to achieve global exponential tracking of the heat dynamics. Using boundary state feedback only, a PI control, coupled to a discontinuous integral term, is designed. Such a controller is capable of  compensating model uncertainties and Lipschitz-in-time disturbances using a continuous control signal. This strategy can be viewed as the merge of the classical integral action and the discontinuous functions commonly used in sliding mode control applications. Using a well-known strategy for motion planning, reference profiles are obtained for the nominal heat equation depending on the choice of the flat output reference. In the presence of parameter uncertainties,  Lyapunov analysis yields the control synthesis, resulting in both norm-bounded tracking errors and  exponential Input-to-State Stability (eISS) of the closed-loop system. Supporting simulations illustrate the closed-loop performance  for different kinds of references in the presence of unbounded but Lipschitz-in-time disturbances and parameter uncertainties.

The outline of this work is as follows. The underlying heat conduction model  and the control objective are introduced in Section \ref{sec:problem}. Using the flat approach, the trajectory generation and open-loop tracking control  are described in Section \ref{sec:gen}. The feedback control design for the nominal and disturbed error systems is given in Section \ref{sec:control}. The reliability of the proposed control strategy is supported in the simulation study of Section \ref{sec:sim}. Finally, some concluding remarks are collected in Section \ref{sec:conclusions}.

\noindent \textbf{Notation:} The function $\lceil\cdot\rfloor^{\gamma}:=|\cdot|^{\gamma}\sign(\cdot)$ is determined for any $\gamma\in \mathbb{R}_{\geq 0}$. Given a differentiable function $r(t)$, notation $r^{(i)}(t)$ with $i \in \mathbb{Z}_{\geq 0}$ stands for the i-th time derivative of $r(t)$.  The Sobolev space of absolutely continuous scalar functions $u(x)$ on $(a, b)$ with square integrable derivatives $u^{(i)}(x),\ i=1,l\ldots,l$ and the norm $\norm{u(\cdot)}_{H^l(a,b)} = \sqrt{\int_a^b \sum_{i=0}^l [u^{(i)}(x)]^2\, dx}$ is typically denoted by $H^l(a,b)$ with $a \leq b$ and $l = \left\{0, 1, 2,...\right\}$. For ease of reference, the nomenclature  $\norm{u(\cdot)}_{H^0(a,b)}=\norm{u(\cdot)}_{L_2(a,b)}=\norm{u(\cdot)}$ is used  throughout. The spatial derivatives are denoted by $u_x =\nicefrac{\partial u}{\partial x}$ and $u_{xx} =\nicefrac{\partial^2 u}{\partial x^2}$.

%For later use, well-known inequalities are recalled.

%\noindent \textbf{Young's Inequality:}
%\begin{equation*}
%  ab \leq \frac{a^p}{p}+\frac{b^q}{q}, \quad p,q>1, \quad \frac{1}{p}+\frac{1}{q}=1.
%\end{equation*}

%\noindent \textbf{Cauchy-Schwarz Inequality:}
%\begin{equation*}
%  \int_a^b f(x)g(x) \, dx \leq \norm{f(x)}\norm{g(x)}.
%\end{equation*}

%\noindent \textbf{Poincare's Inequality:}
%Let $u(x)\in H^1(0,1)$. Then , the following inequality holds
%\begin{equation*}
%  \int_0^1 u^2(x) \, dx \leq 2u^2(i) + 2 \int_0^1 u_x^2(x) \, dx , \quad i=0,1.
%\end{equation*}
%For $u(x)\in H^1(0,D)$ and $D \in \mathbb{R}_{> 0}$, the latter inequality is specified to
%\begin{equation*}
%  \int_0^D u^2(x) \, dx \leq 2Du^2(i) + 2D^2 \int_0^D u_x^2(x) \, dx , \quad i=0,D.
%\end{equation*}

\section{PROBLEM STATEMENT}
\label{sec:problem}

Consider the following heat equation with Neumann boundary conditions (BC)
\begin{equation}
\begin{split}
  u_{t}(x,t) &= d u_{xx}(x,t), \\
  u_x(0,t) &= 0, \quad
  u_x(D,t) = q+\varphi(t),
\end{split}
\label{eq:diff}
\end{equation}
where $u(x,t)$ is the state vector evolving in the space $H^1(0,D)$, $x\in [0,D]$ is the space variable, $t \geq 0$ is the time variable, $d$ is the thermal diffusivity, $q$ is the boundary control input and the function $\varphi(t)$ is a disturbance term supposed to satisfy
\begin{equation}
  \abs{\dot{\varphi}} \leq L
\label{eq:pbound}
\end{equation}
with an \textit{a priori} known constant $L$. Furthermore, the positive parameters $d,D$ are uncertain, but bounded as
\begin{equation}
  0<d_{m} \leq d \leq d_{M}, \quad 0<D_{m} \leq D \leq D_{M},
\label{eq:ubound}
\end{equation}
by some known constant $d_{m},d_{M},D_{m},D_{M}$. Apart from these bounds, The nominal  values $d_n,D_n$ of $d,D$, satisfy  the same bounds  \eqref{eq:ubound}.

It is well-known (see, \textit{e.g.}, \cite[Chapter 3]{b:Orlov-2020}) that for arbitrary initial conditions $u(x,0)$ of class $H^1(0,D)$, there exists a mild solution of the open-loop boundary value problem (BVP) \eqref{eq:diff}. Throughout the paper, only mild solutions of the corresponding BVP are in play.

The objective of this work is to design a control input $q$ capable of driving the output
\begin{equation}
  y(t)=u(0,t)
  \label{eq:output}
\end{equation}
of the underlying  BVP \eqref{eq:diff} to a desired reference $r(t)$, despite the presence of uncertainties and/or disturbances.

\section{MOTION PLANNING}
\label{sec:gen}

Following  \cite{b:Laroche-Martin-Rouchon-2000},\cite[Chapter 12]{b:Krstic-Smyshlyaev-2008}, the reference trajectory generation for the heat equation \eqref{eq:diff} becomes available through its flat output \eqref{eq:output}. The trajectory generation is further performed for the unperturbed system with $\varphi(t)\equiv 0$. To begin with, the perfect knowledge of $d,D$ is assumed.

The state trajectory to follow is represented in the form $\bar{u}(x,t) = \sum_{i=0}^{\infty}a_i(t)\frac{x^i}{i!}$, where the time-varying coefficients $a_i(t)$ have to be determined by substituting the latter sum into \eqref{eq:diff},\eqref{eq:output} and using the desired tracking $r(t)=\bar{u}(0,t)$. Thus, the reference state trajectory is specified to
\begin{equation}
  \bar{u}(x,t) = \sum_{i=0}^{\infty}r^{(i)}(t)\frac{x^{2i}}{d^i(2i)!},
  \label{eq:ref}
\end{equation}
and the nominal input signal from the BC at $x=D$ is defined as
\begin{equation}
  \bar{q} = \sum_{i=1}^{\infty}r^{(i)}(t)\frac{D^{2i-1}}{d^i(2i-1)!}.
  \label{eq:qn}
\end{equation}

In order to guarantee that $u(x,t) \rightarrow \bar{u}(x,t)$ and $y(t) \rightarrow r(t)$, the convergence of \eqref{eq:ref} is to be guaranteed. The next theorem states which conditions should be imposed on the to-be-tracked reference in order to achieve this.

\begin{theorem}\label{th1}
The series \eqref{eq:ref} is absolutely convergent provided that the output reference $r(t)$ fulfils
\begin{equation}
  \sup_{t \geq 0}\abs{r^{(i+1)}(t)} \leq \frac{2d}{D^2}(i+1) \sup_{t \geq 0}\abs{r^{(i)}(t)} \quad \forall \quad i \in \mathbb{Z}_{\geq 0}.
  \label{eq:refselect}
\end{equation}
\end{theorem}

\begin{proof}
In order to check if under conditions of the theorem, \eqref{eq:ref} is a convergent series, the ratio test is performed with $S_p = r^{(p)}(t)\frac{x^{2p}}{d^p(2p)!}$. Then
\begin{equation*}
\begin{split}
  \lim_{p\to\infty} \abs{\frac{S_{p+1}}{S_p}} &= \lim_{p\to\infty} \abs{\frac{r^{(p+1)}(t)}{r^{(p)}(t)}} \frac{x^{2p+2}}{x^{2p}} \frac{d^{p}}{d^{p+1}} \frac{(2p)!}{(2p+2)!}\\
  &= \lim_{p\to\infty} \frac{x^2}{d} \abs{\frac{r^{(p+1)}(t)}{r^{(p)}(t)}} \frac{1}{(2p+1)(2p+2)}<1
\end{split}
\end{equation*}
for all $x \in [0,D]$ under condition \eqref{eq:refselect}. Theorem \ref{th1} is thus proved.
\end{proof}

References $r(t)$, which are usually adopted in motion planning \cite{b:Laroche-Martin-Rouchon-2000,b:Lynch-Rudolph-2002,b:Krstic-Magnis-Vazquez-2009,b:Meurer-Kugi-2009}, are the Gevrey class defined as follows.

\begin{definition}\cite{b:Laroche-Martin-Rouchon-2000}
A smooth function $r(t)$ is Gevrey of order $\alpha$ if exist $M,R >0$ such that $\sup_{t \geq 0} \abs{r^{(i)}(t)} \leq M \frac{i!^\alpha}{R^i}$, for all $i \in \mathbb{Z}_{\geq 0}$.
\label{def:Gevrey}
\end{definition}

In the works cited above, the condition to guarantee the series convergence requires  $r(t)$ to be Gevrey of order $\alpha<2$. The next lemma links references $r(t)$ satisfying  condition  \eqref{eq:refselect} to a certain kind of the Gevrey  functions.

\begin{lemm}
Any smooth function $r(t)$, fulfilling \eqref{eq:refselect}, is Gevrey of order $\alpha=1$ such that Definition \ref{def:Gevrey} holds for $r(t)$ with $R=\frac{D^2}{2d}$ and $M>0$.
\end{lemm}

\begin{proof}
From the Definition \ref{def:Gevrey}, one concludes that
\begin{equation*}
\begin{split}
  \sup_{t \geq 0} \abs{r^{(i+1)}(t)} \leq M \frac{(i+1)!^\alpha}{R^{i+1}}
  \leq \frac{M}{R}(i+1)^{\alpha} \frac{i!^\alpha}{R^{i}},
\end{split}
\end{equation*}
thereby deriving the inequality
\begin{equation*}
  \sup_{t \geq 0} \abs{r^{(i+1)}(t)} \leq \frac{(i+1)^{\alpha}}{R}\sup_{t \geq 0} \abs{r^{(i)}(t)}.
\end{equation*}
Since  \eqref{eq:refselect} follows from the latter inequality with the selection of $\alpha,R,M$ under the lemma conditions, the proof is thus completed.
\end{proof}

Although the principle reference, used in the afore-cited works, is the so-called "bump function" (see \cite{b:Laroche-Martin-Rouchon-2000}), its proposed counterpart \eqref{eq:refselect} allows one to exemplify more admissible references among analytical  functions such as:
\begin{enumerate}
\renewcommand{\labelenumi}{\textbf{F.\theenumi}}
  \item $\boldsymbol{r(t)=A}$ for any $A\in \Re$.
  \item $\boldsymbol{r(t)=At}$, $t \in [0,t_f]$ for any $A\in \Re$, $t_f \geq \frac{D^2}{2d}$.
  \item $\boldsymbol{r(t)=A \sin(\omega t)}$, $\omega \in \Re_{>0}$ for any $A\in \Re$, $\omega \leq \frac{2d}{D^2}$.
  \item $\boldsymbol{r(t)=A e^{\beta t}}$, $\beta \in \Re$ for any $A\in \Re$, $\abs{\beta} \leq \frac{2d}{D^2}$.
\end{enumerate}

\begin{remark}
The output \eqref{eq:output} matches the reference  trajectory $r(t)$ only if  the state trajectory to follow satisfies the same initial condition as that of the plant, \textit{i.e.}, if $u(x,0)=\bar{u}(x,0)$. Furthermore, the nominal control \eqref{eq:qn} is not capable of  compensating any non-trivial disturbance  $\varphi(t) \neq 0$. Next  the tracking problem of interest is addressed for the arbitrarily initialized state trajectory to follow in the presence of matched disturbances.
\end{remark}

\section{TRACKING CONTROL DESIGN}
\label{sec:control}

Setting the state deviation
\begin{equation}
  \tilde{u}(x,t) = u(x,t)-\bar{u}(x,t),
  \label{eq:error}
\end{equation}
from the reference trajectory $\bar{u}(x,t)$, given by \eqref{eq:ref}, the error dynamics \eqref{eq:error} are then governed by
\begin{equation}
\begin{split}
  \tilde{u}_{t}(x,t) &= d \tilde{u}_{xx}(x,t), \\
  \tilde{u}_x(0,t) &= 0, \quad 
  \tilde{u}_x(D,t) = q-\bar{q}+\varphi(t)=\tilde{q}+\varphi(t),
\end{split}
\label{eq:diffe}
\end{equation}
where the control input
\begin{equation}
  q = \tilde{q}+\bar{q}
  \label{eq:q}
\end{equation}
is precomposed in terms of the reference input $\bar{q}$, determined by \eqref{eq:qn}, and the virtual component $\tilde{q}$.
Now the tracking objective  for the nominal reference \eqref{eq:ref} is reduced  to the virtual input design $\tilde{q}$ exponentially stabilizing the error dynamics  \eqref{eq:diffe} in the origin.

\subsection{Disturbance-free tracking}

First, the disturbance-free case  $\varphi(t)\equiv 0$ is analysed. Selecting the control $\tilde{q}$ as
\begin{equation}
  \tilde{q} = - \frac{\lambda_1}{D_n} \tilde{u}(D,t),
  \label{eq:qt1}
\end{equation}
where $\lambda_1$ is a gain to be tuned and $D_n$ the nominal value of $D$. The next result is in order.

\begin{theorem}
Let  $\varphi(t)\equiv 0$ and let the control parameters in \eqref{eq:qt1} be such that 
\begin{equation}
  \lambda_1 > \frac{D_n}{D_{m}}.
  \label{eq:lambda1}
\end{equation}
Then the error dynamics \eqref{eq:diffe}, driven by  \eqref{eq:qt1}, are globally exponentially stable.
\hfill $\blacksquare$
\end{theorem}

\begin{proof}
Consider the  positive definite Lyapunov functional candidate $V = \frac{1}{2}\norm{\tilde{u}(x,t)}^2$. Its derivative along the system \eqref{eq:diffe} reads as
\begin{equation*}
\begin{split}
  \dot{V} = \int_0^D \tilde{u}(x,t)\tilde{u}_t(x,t) \, dx
  = d \int_0^D \tilde{u}(x,t)\tilde{u}_{xx}(x,t) \, dx .
\end{split}
\end{equation*}
Applying the integration  by parts and then  employing the BC and  Poincare's inequality, it follows that 
\begin{equation*}
\begin{split}
  \dot{V} &= -d \int_0^D \tilde{u}_{x}^2(x,t) \, dx + d \left[\tilde{u}(x,t)\tilde{u}_x(x,t)\right]_0^D \\
  &\leq -\frac{d}{2 D^2} \norm{\tilde{u}(x,t)}^2 + \frac{d}{D}\tilde{u}^2(D,t) + d\tilde{u}(D,t)\tilde{q}.
\end{split}
\end{equation*}

Now, taking into account the control law \eqref{eq:qt1}, coupled to \eqref{eq:lambda1}, the Lyapunov derivative is further estimated as
\begin{equation*}
\begin{split}
  \dot{V} &\leq -\frac{d}{2 D^2} \norm{\tilde{u}(x,t)}^2 - \left(\frac{\lambda_1}{D_n}-\frac{1}{D}\right)d\tilde{u}^2(D,t)\\
  &\leq -\frac{d}{2 D^2} \norm{\tilde{u}(x,t)}^2
  \leq -\frac{d}{D^2} V,
\end{split}
\end{equation*}
that guarantees the exponential decay of the system dynamics \eqref{eq:diffe}.  Since the Lyapunov functional is radially unbounded, \eqref{eq:diffe} is concluded to be globally exponentially stable.
\end{proof}

\subsection{Tracking under disturbances}

In the presence of the external disturbance $\varphi(t)\neq 0$, the proposed control $\tilde{q}$ in \eqref{eq:qt1} is no longer capable of  stabilizing the error dynamics. To robustify the control law in the disturbance-corrupted case, it is modified to
\begin{equation}
\begin{split}
  \tilde{q} = - \frac{\lambda_1}{D_n} \tilde{u}(D,t) + \nu, \quad
  \dot{\nu} = - \lambda_2 \tilde{u}(D,t) - \lambda_3 \Sabs{\tilde{u}(D,t)}^0,
\end{split}
  \label{eq:qt2}
\end{equation}
where $\lambda_1,\lambda_2,\lambda_3$ are gains to be tuned and $D_n$ is the nominal value of $D$.

The proposed feedback law \eqref{eq:qt2} is composed by a PI control and a discontinuous term passing through an integrator. It generates a continuous control signal despite having a discontinuous (multi-valued) right-hand side in the manifold $\tilde{u}(D,\cdot)=0$. The precise meaning of the solutions of the distributed parameter system \eqref{eq:diffe} driven by this discontinuous controller, are viewed in the Filippov sense \cite{b:filippov}. Extension of the Filippov concept towards the infinite-dimensional setting may be found in \cite{b:Orlov-2020}. The present paper focuses on the tracking synthesis whereas the well-posedness analysis of the closed-loop system \eqref{eq:diffe},\eqref{eq:qt2} is similar to that of \cite{b:Pisano-Orlov-Usai-2011} and it remains beyond the scope of the paper. Thus, for the closed-loop system in question, it is assumed that it possesses a  Filippov solution.

The closed-loop system \eqref{eq:diffe}, driven by \eqref{eq:qt2}, reads as
\begin{equation}
\begin{split}
  \tilde{u}_{t}(x,t) &= d \tilde{u}_{xx}(x,t), \\
  \tilde{u}_x(0,t) &= 0, \quad
  \tilde{u}_x(D,t) = - \frac{\lambda_1}{D_n} \tilde{u}(D,t) + \delta, \\
  \dot{\delta} &= - \lambda_2 \tilde{u}(D,t) - \lambda_3 \Sabs{\tilde{u}(D,t)}^0 + \dot{\varphi}(t),
\end{split}
\label{eq:diffe2}
\end{equation}
with $\delta = \nu + \varphi$, being substituted into \eqref{eq:qt2} for $\nu= \delta - \varphi$ for deriving $\delta$-dynamics \eqref{eq:diffe2}.  The following result  is then in force.

\begin{theorem}
 Let the error dynamics \eqref{eq:diffe} be driven by the feedback  $\tilde{q}$, governed by \eqref{eq:qt2} and tuned in accordance with
\begin{equation}
  \lambda_1 > \frac{D_n}{D_{m}}, \quad \lambda_2 > d_{M}, \quad \lambda_3 > L.
  \label{eq:lambda2}
\end{equation}
Then the  closed-loop error dynamics \eqref{eq:diffe2} are globally exponentially stable despite the presence of any boundary  disturbance $\varphi$ of class \eqref{eq:pbound}.
\end{theorem}

\begin{proof}
Consider  the  positive definite Lyapunov functional candidate
\begin{equation}
  V = \frac{1}{2}\norm{\tilde{u}(x,t)}^2 + \frac{1}{2}\delta^2,
  \label{eq:Lyap}
\end{equation}
and employing the magnitude bound \eqref{eq:pbound} for $\dot{\varphi}$, compute  its time derivative along the error dynamics  \eqref{eq:diffe2}, thus obtaining
\begin{equation*}
\begin{split}
  \dot{V} &\leq -\frac{d}{2 D^2} \norm{\tilde{u}(x,t)}^2 - \left(\frac{\lambda_1}{D_n}-\frac{1}{D}\right)d\tilde{u}^2(D,t)\\
  &\quad + d\tilde{u}(D,t)\delta + (- \lambda_2 \tilde{u}(D,t) - \lambda_3 \Sabs{\tilde{u}(D,t)}^0 + \dot{\varphi})\delta\\
  &\leq -\frac{d}{2 D^2} \norm{\tilde{u}(x,t)}^2 - \left(\frac{\lambda_1}{D_n}-\frac{1}{D}\right)d\tilde{u}^2(D,t)\\
  &\quad + (- \lambda_2 +d - \lambda_3 \abs{\tilde{u}(D,t)}^{-1} + L\Sabs{\tilde{u}(D,t)}^{-1})\tilde{u}(D,t)\delta.
\end{split}
\end{equation*}
 By applying Young's inequality and the gains selection  \eqref{eq:lambda2}, it follows that
\begin{equation*}
\begin{split}
  \dot{V} &\leq -\frac{d}{2 D^2} \norm{\tilde{u}(x,t)}^2 - \left(\frac{\lambda_1}{D_n}-\frac{1}{D}\right)d\tilde{u}^2(D,t)\\
  &\quad - (\lambda_2 -d + \lambda_3 \abs{\tilde{u}(D,t)}^{-1} - L\Sabs{\tilde{u}(D,t)}^{-1})\frac{1}{2}\delta^2\\
  &\quad - (\lambda_2 -d + \lambda_3 \abs{\tilde{u}(D,t)}^{-1} - L\Sabs{\tilde{u}(D,t)}^{-1})\frac{1}{2}\tilde{u}^2(D,t)\\
  &\leq -\frac{d}{2 D^2} \norm{\tilde{u}(x,t)}^2\\
  &\quad - (\lambda_2 -d + \lambda_3 \abs{\tilde{u}(D,t)}^{-1} - L\Sabs{\tilde{u}(D,t)}^{-1})\frac{1}{2}\delta^2\\
  &\leq - \alpha V
\end{split}
\end{equation*}
where $\alpha = \mathrm{min}\left\{\frac{d}{D^2},\lambda_2-d+\lambda_3-L \right\}$. The latter inequality ensures the exponential stability of  \eqref{eq:diffe2}. Since the Lyapunov functional is radially unbounded, the result holds globally.
\end{proof}

\begin{remark}
Due to the exponential decay of the Lyapunov functional \eqref{eq:Lyap} and by virtue of $\delta=\nu+\varphi$, the virtual input $\nu$ approaches the negative disturbance value, \textit{i.e.}, $\nu(t)\rightarrow -\varphi(t)$ as $t \rightarrow \infty$.
\end{remark}

\begin{remark}
The proposed feedback \eqref{eq:qt2} requires only the boundary state information $u(D,t)$ at $x=D$ only. Thus, the boundary output feedback is available to perform the tracking task.
\end{remark}

\subsection{Tracking under uncertain motion planning}

The feedback controller, constructed by now,  enforces  the heat conduction dynamics to exponentially track an output reference signal. An admissible reference to be tracked  is obtained from the motion planning assuming to perfectly know the nominal system \eqref{eq:diff}. This is not however the case in realistic applications where only nominal values $d_n$, $D_n$ of respectively $d$, $D$ are available. From now on, the reference state trajectory $\bar{u}(x,t) = \sum_{i=0}^{\infty}r^{(i)}(t)\frac{x^{2i}}{d_n^i(2i)!}$, and the nominal input signal $\bar{q} = \sum_{i=1}^{\infty}r^{(i)}(t)\frac{D_n^{2i-1}}{d_n^i(2i-1)!}$, rely on the nominal plant parameter values rather than on their real values.

Using the same error variable  \eqref{eq:error}, the closed-loop  error dynamics, enforced by \eqref{eq:qt2},  are governed by
\begin{equation}
\begin{split}
  \tilde{u}_{t}(x,t) &= d \tilde{u}_{xx}(x,t) + M(x,t), \\
  \tilde{u}_x(0,t) &= 0, \quad
  \tilde{u}_x(D,t) = - \frac{\lambda_1}{D_n} \tilde{u}(D,t) + \delta, \\
  \dot{\delta} &= - \lambda_2 \tilde{u}(D,t) - \lambda_3 \Sabs{\tilde{u}(D,t)}^0 + \dot{\varphi}(t),
\end{split}
\label{eq:diffe3}
\end{equation}
where the perturbation term 
\begin{equation}
\begin{split}
  M(x,t) = (d-d_n)\sum_{i=1}^{\infty}r^{(i)}(t)\frac{x^{2i-2}}{d_n^i(2i-2)!}
\end{split}
\end{equation}
is due to the uncertainties on $d$. 
\begin{lemm}
The uncertain term $M(x,t)$ is uniformly bounded 
\begin{equation}
  \norm{M(\cdot,t)} \leq (d_{M}-d_n)^2 L_1 \quad \forall \quad t\geq 0,
  \label{eq:Mbound}
\end{equation}
with some $L_1 \in \Re_{\geq 0}$ provided that the output reference $r(t)$ fulfils the condition
\begin{equation}
  \sup_{t \geq 0} \abs{r^{(i+1)}(t)} \leq \frac{d_n}{D_M^2} \sup_{t \geq 0} \abs{r^{(i)}(t)}, \quad \forall \quad i \in \mathbb{Z}_{\geq 0}.
  \label{eq:errorcond}
\end{equation}
\end{lemm}

\begin{proof}
With the upper bound $D_M$ of $d$, the norm of $M(x,t)$ is estimated as follows
\begin{equation*}
\begin{split}
  \norm{M(x,t)} &= \int_0^D M^2(x,t) \, dx \\
  &\leq (d_{M}-d_n)^2 \sum_{i=1}^{\infty}[r^{(i)}(t)]^2 \int_0^D \frac{x^{4i-4}}{d_n^{2i}(2i-2)!^2} \, dx \\
  &= (d_{M}-d_n)^2 \sum_{i=1}^{\infty}[r^{(i)}(t)]^2 \frac{D^{4i-3}}{d_n^{2i}(4i-3)(2i-2)!^2}.
\end{split}
\end{equation*}

In order to check if the latter series is absolutely convergent, the ratio test is now performed  with $S_p=[r^{(p)}(t)]^2 \frac{D^{4p-3}}{d_n^{2p}(4p-3)(2p-2)!^2}$ to conclude that
{\footnotesize
\begin{equation*}
\begin{split}
  \lim_{p\to\infty} \abs{\frac{S_{p+1}}{S_p}} &= \lim_{p\to\infty} \left[ \frac{r^{(p+1)}(t)}{r^{(p)}(t)}\right]^2 \frac{D^{4p+1}}{D^{4p-3}} \frac{d_n^{2p}}{d_n^{2p+2}} \abs{\frac{4p-3}{4p+1}} \left[ \frac{(2p-2)!}{(2p)!}\right]^2 \\
  &= \lim_{p\to\infty} \left[ \frac{D^2}{d_n} \frac{r^{(p+1)}(t)}{r^{(p)}(t)}\right]^2 \left[\frac{1}{(2p-1)(2p)}\right]^2 \abs{\frac{4p-3}{4p+1}}<1
\end{split}
\end{equation*}}
provided condition \eqref{eq:errorcond} holds true. Hence, the series is convergent and approaches a finite value $L_1$, thereby ensuring  the term $\norm{M(x,t)}$ is bounded according to \eqref{eq:Mbound}.
\end{proof}

For later use, let us specify the eISS concept for solutions $[\tilde{u}(x,t),\delta]$ of the BVP \eqref{eq:diffe3} with respect to the $H^0(0,D)$-norm
\begin{equation}
  \Gamma(t) = \left(\norm{\tilde{u}(x,t)}^2+ \delta^2 \right)^{1/2};
  \label{eq:norm}
\end{equation}
see \cite{b:Dashkovskiy-Mironchenko-2013} for ISS details  in the infinite-dimensional setting.
\begin{definition}
System \eqref{eq:diffe3} is called eISS in the sense of the norm \eqref{eq:norm} with respect to $M(x,t)$ if there exist $\alpha_1,\alpha_2>0$ and a function $\gamma$ of class $\mathscr{K}$ such that the inequality
\begin{equation}
  \Gamma(t) \leq \alpha_1 e^{-\alpha_2 t} \Gamma(0) + \gamma(L_1),
  \label{eq:ISS}
\end{equation}
holds on arbitrary solution $\tilde{u}(x,t)$ of \eqref{eq:diffe3}, with $\Gamma(0)< \infty$ and $M(x,t):\norm{M(\cdot,t)}\leq L_1 < \infty$ for all $t\geq 0$.
\end{definition}

The eISS  of system \eqref{eq:diffe3} is then established.
\begin{theorem}
 Let the controller gains in \eqref{eq:qt2} be such that \eqref{eq:lambda2} holds and let the output reference $r(t)$ respect  condition \eqref{eq:Mbound}. Then system \eqref{eq:diffe3} is eISS. 
\end{theorem}

\begin{proof}
The norm of the Lyapunov functional \eqref{eq:Lyap} is straightforwardly estimated as $m_1 \Gamma^2(t) \leq V(t) \leq m_2 \Gamma^2(t)$, with some $m_1,m_2>0$. Following the same steps as before, its derivative along the trajectories of system \eqref{eq:diffe3} is evaluated as
\begin{equation*}
\begin{split}
  \dot{V} &\leq -\frac{d}{2 D^2} \norm{\tilde{u}(x,t)}^2 + \int_0^D \tilde{u}(x,t)M(x,t) \, dx\\
  &\quad - (\lambda_2 -d + \lambda_3 \abs{\tilde{u}(D,t)}^{-1} - L\Sabs{\tilde{u}(D,t)}^{-1})\frac{1}{2}\delta^2.
\end{split}
\end{equation*}

Using the Cauchy-Schwarz inequality and the above bounds of the Lyapunov functional as well as bound \eqref{eq:Mbound} of the term $\norm{M(x,t)}$, its derivative reads as
\begin{equation*}
\begin{split}
  \dot{V} &\leq -\alpha m_2 \Gamma^2(t) + \norm{\tilde{u}(x,t)}\norm{M(x,t)}\\
  &\leq -\alpha m_2 (1-\theta) \Gamma^2(t) -\alpha m_2 \theta \Gamma^2(t) + (d_{M}-d_n)^2 L_1 \Gamma(t) \\
  &\leq -\alpha m_2 (1-\theta) \Gamma^2(t) \quad \forall \quad \Gamma(t) \geq \frac{(d_{M}-d_n)^2}{\alpha m_2 \theta}L_1,
\end{split}
\end{equation*}
where $0<\theta<1$, $\alpha = \mathrm{min}\left\{\frac{d}{D^2},\lambda_2-d+\lambda_3-L \right\}$. Last expression results in the norm $\Gamma(t)$ to be globally ultimately bounded with respect to $M(x,t)$. Furthermore, applying the comparison lemma (see \cite{b:Khalil2002}) one concludes the  eISS condition \eqref{eq:ISS} with $\gamma(L_1)=\frac{(d_{M}-d_n)^2}{\alpha m_2 \theta}L_1$. This completes the proof. 
\end{proof}

\section{SIMULATIONS}
\label{sec:sim}

The proposed control strategy \eqref{eq:qn},\eqref{eq:q},\eqref{eq:qt2} has been implemented in the system \eqref{eq:diff} using Matlab Simulink with Euler's integration method of fixed step and a sampling time equal to $50$ [ms]. The heat equation was implemented using $d=0.05$, $D=10$ and $\varphi(t)=0.01t+2\sin(t)$ and the finite-differences approximation technique, discretizing the spatial domain $x\in[0,D]$ into 51 ordinary differential equations (ODE). The initial condition was set in $u(x,0)=0$.

The tracking reference $r(t)$ and the corresponding nominal reference $\bar{u}(x,t)$ and control $\bar{q}$, are defined depending on the simulation time and expressions \eqref{eq:ref}-\eqref{eq:qn}:
\begin{enumerate}
  \item \textbf{$\boldsymbol{t \in [0,1\times 10^4)}$ [s]:} \\
  $r(t)=A t$, $A=1\times 10^{-3}$, \\
  $\bar{u}(x,t)=At + A \frac{x^2}{2d}$, 
  $\bar{q} = A \frac{D}{d}$.
  \item \textbf{$\boldsymbol{t \in [1\times 10^4,2\times 10^4)}$ [s]:} \\
  $r(t)=20$, 
  $\bar{u}(x,t)= 20$, 
  $\bar{q} = 0$.
  \item \textbf{$\boldsymbol{t \in [2\times 10^4,3\times 10^4)}$ [s]:} \\
  $r(t)=A e^{[-\beta (t-2\times 10^4)]}$, $A=20$, $\beta=1\times 10^{-3}$\\
  $\bar{u}(x,t)= A e^{[-\beta (t-2\times 10^4)]}\cos\left(x\sqrt{\frac{\beta}{d}} \right)$, \\
  $\bar{q} = -A \sqrt{\frac{\beta}{d}} e^{[-\beta (t-2\times 10^4)]}\sin\left(D\sqrt{\frac{\beta}{d}} \right)$.
  \item \textbf{$\boldsymbol{t \in [3\times 10^4,4\times 10^4]}$ [s]:} \\
  $r(t)=A \sin(\omega t)$, $A=10$, $\omega=1\times 10^{-3}$, \\
  $\bar{u}(x,t)=\frac{1}{2}A e^{x\sqrt{\frac{\omega}{2d}}}\sin \left(\omega t + x\sqrt{\frac{\omega}{2d}} \right)$ \\ \hspace*{37pt}$+\frac{1}{2}A e^{-x\sqrt{\frac{\omega}{2d}}}\sin \left(\omega t - x\sqrt{\frac{\omega}{2d}} \right)$, \\
  $\bar{q} = \frac{1}{2}A \sqrt{\frac{\omega}{2d}} e^{D\sqrt{\frac{\omega}{2d}}}\sin \left(\omega t + D\sqrt{\frac{\omega}{2d}} \right)$ \\
  \hspace*{17pt}$+\frac{1}{2}A \sqrt{\frac{\omega}{2d}} e^{D\sqrt{\frac{\omega}{2d}}}\cos \left(\omega t + D\sqrt{\frac{\omega}{2d}} \right)$ \\
  \hspace*{16pt}$-\frac{1}{2}A \sqrt{\frac{\omega}{2d}} e^{-D\sqrt{\frac{\omega}{2d}}}\sin \left(\omega t - D\sqrt{\frac{\omega}{2d}} \right)$ \\
  \hspace*{16pt}$-\frac{1}{2}A \sqrt{\frac{\omega}{2d}} e^{-D\sqrt{\frac{\omega}{2d}}}\cos \left(\omega t - D\sqrt{\frac{\omega}{2d}} \right)$.\footnote{See \cite[Chapter 12]{b:Krstic-Smyshlyaev-2008} for more details on how to obtain the analytic nominal expressions for exponential and sinusoidal references.}
\end{enumerate}

\begin{figure}[ht!]
  \centering
  \includegraphics[width=7.8cm,height=3.5cm]{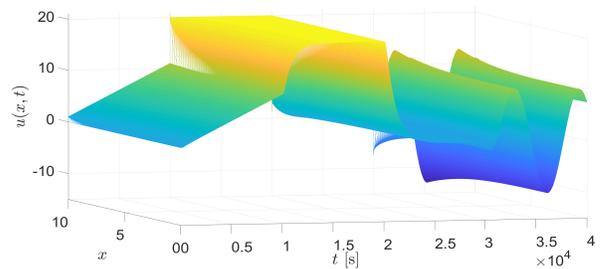}
  \caption{Evolution of the heat equation state $u(t,x)$.}
  \label{fig:u}
\end{figure}

\begin{figure}[ht!]
  %\centering
  \hspace{20pt}
  \includegraphics[width=7.3cm,height=3.5cm]{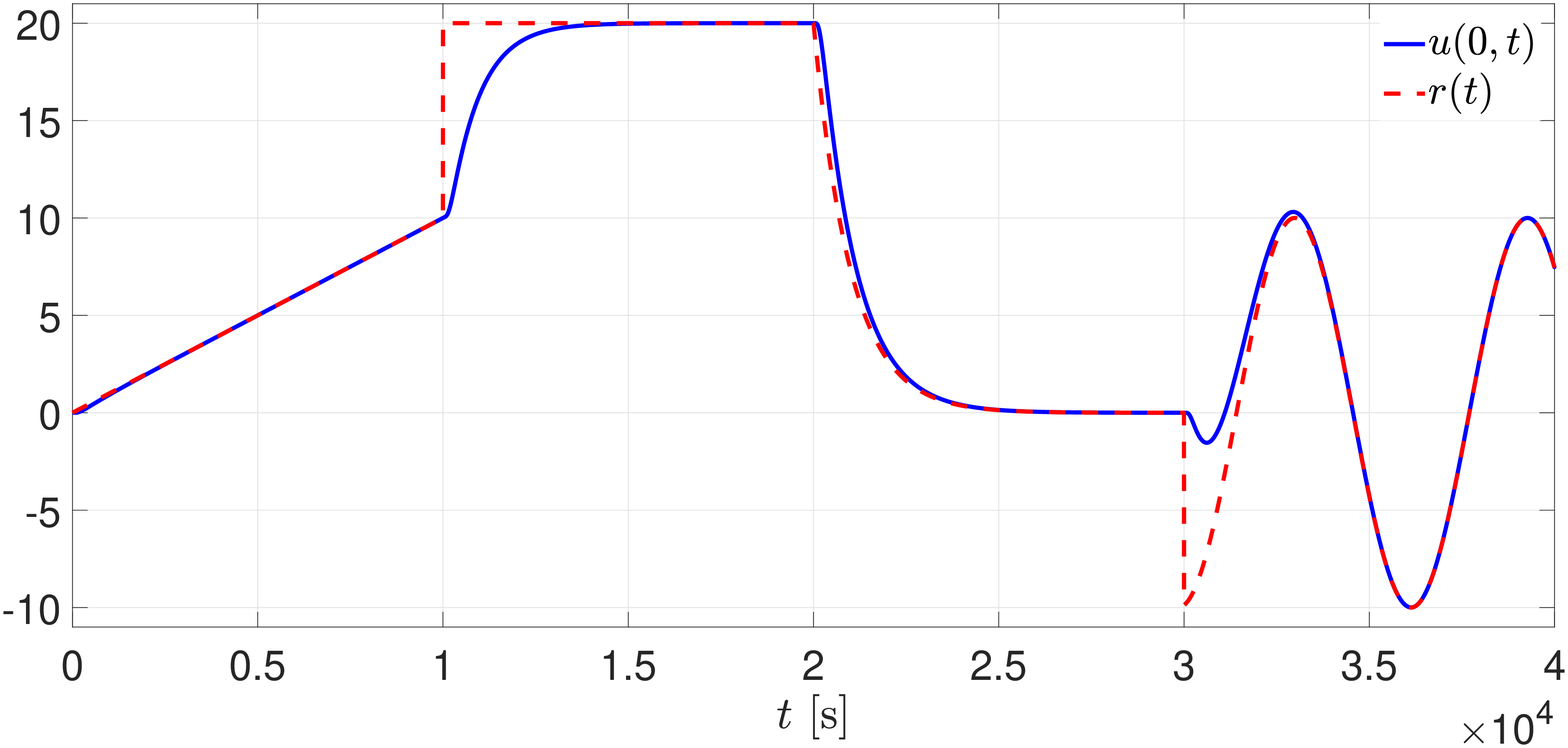}
  \caption{Tracking of the output $y=u(0,t)$ over the successive references $r(t)$.}
  \label{fig:u0}
\end{figure}

\begin{figure}[ht!]
  \centering
  \includegraphics[width=7.8cm,height=3.5cm]{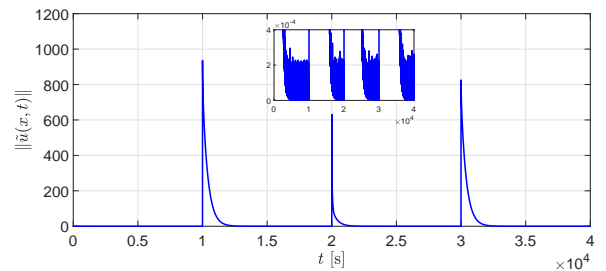}
  \caption{Norm of the error.}
  \label{fig:norm}
\end{figure}

\begin{figure}[ht!]
  \centering
  \includegraphics[width=7.8cm,height=3.5cm]{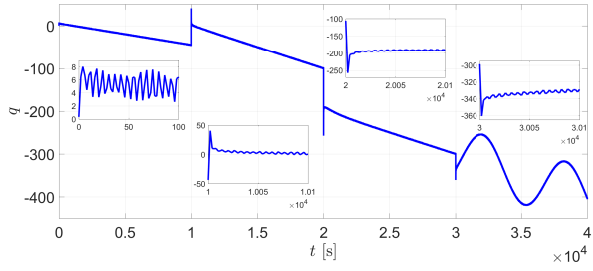}
  \caption{Control signal $q$.}
  \label{fig:q}
\end{figure}

All of the proposed output references fulfil the condition \eqref{eq:refselect}, and more precisely, the conditions described in \textbf{F.1}-\textbf{F.4}. The respective gains of the error control $\tilde{q}$ have been selected according to \eqref{eq:lambda2} as $\lambda_1=1,\lambda_2=10,\lambda_3=2.5$ and $D_n=9$. The results are displayed in Figs. \ref{fig:u}-\ref{fig:q}. The solution $u(x,t)$ of the heat equation is performing a tracking over the space variable $x$. The main objective of stabilizing the output $y=u(0,t)$ over the four successive references $r(t)$ is achieved despite the presence of the unbounded but Lipschitz perturbation $\varphi(t)$ (see Fig. \ref{fig:u0}). Furthermore, the designed control strategy is able to stabilize the norm of the error despite the abrupt change between references, as seen in Fig. \ref{fig:norm}. The boundary control $q$ shown in Fig. \ref{fig:q} clarifies how such disturbance compensation is performed. This is due to the presence of the discontinuous term on the control design. Nevertheless, the control signal generated is continuous throughout the tracking task.

\begin{figure}[ht!]
  %\centering
  \hspace{20pt}
  \includegraphics[width=7.3cm,height=3.5cm]{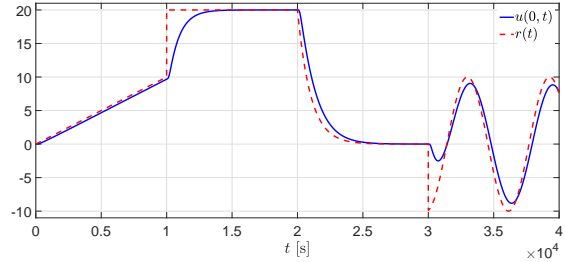}
  \caption{Tracking of the output $y=u(0,t)$ with uncertainties in the trajectory generation.}
  \label{fig:u02}
\end{figure}

\begin{figure}[ht!]
  \centering
  \includegraphics[width=7.8cm,height=3.5cm]{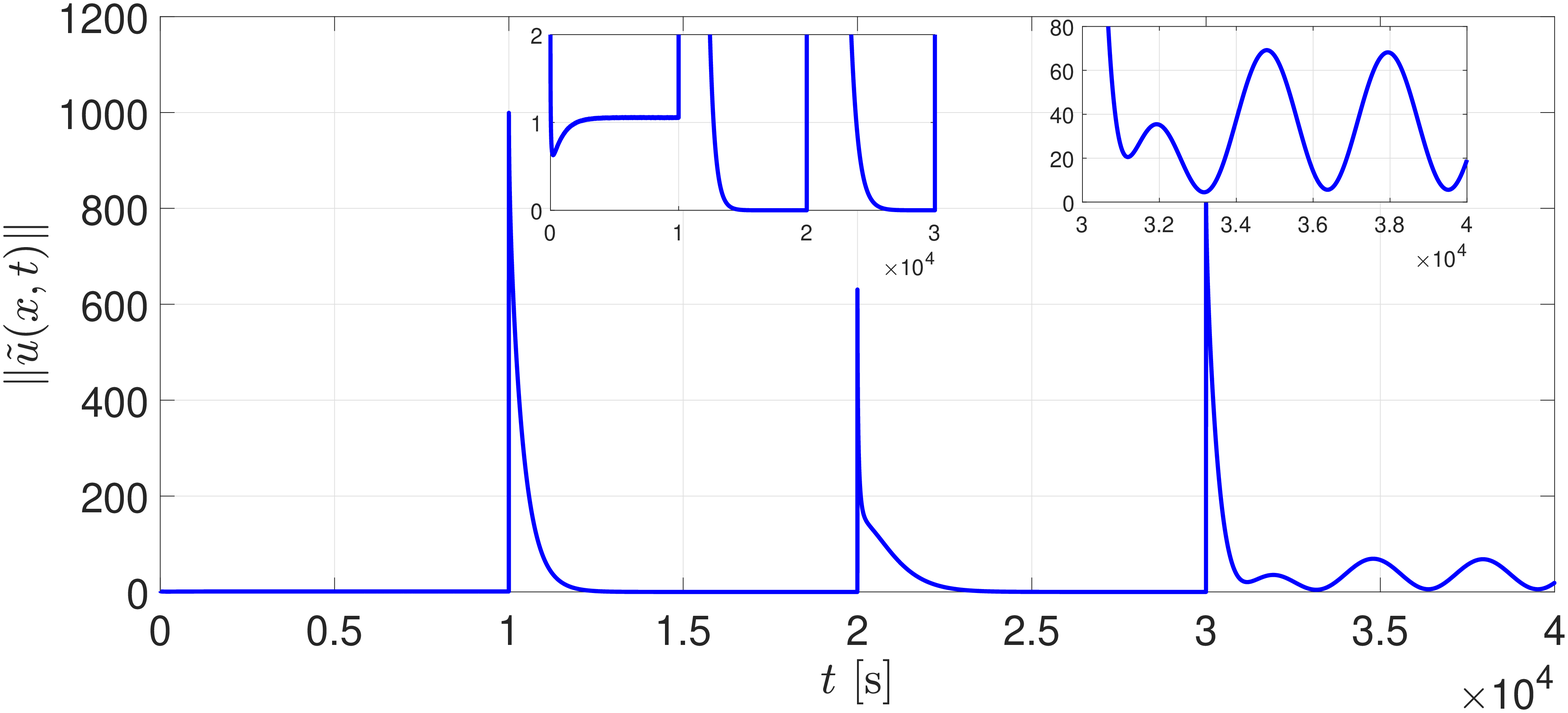}
  \caption{Norm of the error with uncertainties in the trajectory generation.}
  \label{fig:norm2}
\end{figure}

In order to show the scenario where the motion planning was performed used nominal values of $d$ and $D$, the same simulations have been made but using $d_n=0.06$ and $D_n=9$. In this case, the obtained reference $\bar{u}(x,t)$ and nominal control $\bar{q}$ introduce an error. The results are shown in Figs. \ref{fig:u02}-\ref{fig:norm2}. The tracking error is able again to compensate the unbounded disturbance and force the trajectories to follow this new \textit{wrong} reference, that is why the norm of the error is not zero. The magnitude of the error depends on how close the nominal values are from the real system parameters and the kind of reference to be followed, \textit{i.e.}, the bound $\gamma(L_1)$ obtained in the eISS analysis.

\section{CONCLUSIONS}
\label{sec:conclusions}

The heat equation with boundary control is analysed, and robust output tracking  is developed. The proposed boundary control requires the state at the boundary only and it is composed of a PI control and an extra discontinuous term, passing through an integrator. Such a controller is typically used the sliding mode control theory and  it is derived from a Lyapunov approach.  The  controller, thus composed, compensates Lipschitz-in-time disturbances and uncertainties in the system. It  generates sliding modes in the actuator dynamics so that after passing through the integrator, a continuous control signal is applied to the underlying system, thereby diminishing the chattering effect. Capabilities of tracking heat conduction dynamics along  different kinds of boundary references and good robustness properties of the developed design are illustrated in the simulation study.  The reference profiles are obtained for the nominal heat conduction model using the flatness approach.  Being applied to the disturbed plant model with uncertain parameters, the closed-loop  eISS is additionally established along with the boundedness of the tracking error norm. Robustification of the motion planning to inherit the exponential stability from the nominal plant is among open problems calling for further investigation.

                                  % on the last page of the document manually. It shortensplant

%%%%%%%%%%%%%%%%%%%%%%%%%%%%%%%%%%%%%%%%%%%%%%%%%%%%%%%%%%%%%%%%%%%%%%%%%%%%%%%%

%\section*{APPENDIX}

\section*{ACKNOWLEDGEMENT}

The authors would like to acknowledge the support of the European Research Council (ERC) under the European Union's Horizon 2020 research and innovation program (Grant agreement no. 757848 CoQuake). Prof. Yury Orlov work has been partially supported by Atlanstic2020, a research program of R\'egion Pays de la Loire in France, and by CONACYT grant A1-S-9270.

\bibliography{Bibliografias}
\bibliographystyle{IEEEtran}

\end{document}